\documentclass[a4paper]{article}
\usepackage{a4}

\usepackage{amsmath,amssymb,amsthm}
\usepackage{graphicx}
\usepackage{graphics}
\usepackage{epsfig}
\usepackage{hangcaption}
\usepackage{hyperref}

\usepackage{pictexwd}
\usepackage{dcpic}
\usepackage[all]{xy}
\usepackage{amstext}
\usepackage{amsfonts}
\usepackage{latexsym}

\newtheorem{defin}{Definition}[section]
\newtheorem{thm}[defin]{Theorem}
\newtheorem{cor}[defin]{Corollary}
\newtheorem{rem}[defin]{Remark}
\newtheorem{lem}[defin]{Lemma}
\newtheorem{prop}[defin]{Proposition}
\newtheorem{quest}[defin]{Question}
\newtheorem{prob}[defin]{Problem}
\newtheorem{axi}[defin]{Axiom}
\newtheorem{ex}[defin]{Example}
\newtheorem{notice}[defin]{Note}
\newtheorem{conven}[defin]{Convention}
\newtheorem{conj}[defin]{Conjecture}

\newcommand{\theorem}[1]{\begin{thm} \sl{#1} \end{thm}}
\newcommand{\theoremname}[2]{\begin{thm}[#1] \sl{#2} \end{thm}}

\newcommand{\definition}[1]{\begin{defin} \emph{#1} \end{defin}}

\newcommand{\remark}[1]{\begin{rem} \emph{#1} \end{rem}}

\newcommand{\lemmaname}[2]{\begin{lem}[#1] \sl{#2} \end{lem}}

\newcommand{\PL}{\mathrm{PL}}

\DeclareMathOperator{\rot}{\mathrm{rot}}

\DeclareMathSymbol{\Rb}{\mathbin}{AMSb}{"52}
\DeclareMathSymbol{\Zb}{\mathbin}{AMSb}{"5A}

\title{\textsc{Mather invariants in groups of piecewise-linear homeomorphisms}}

\author{
Francesco Matucci\thanks{The author gratefully acknowledges the Centre de Recerca Matem\`atica (CRM)
and its staff for the support received during the development of this work.}
}

\date{}

\begin{document}

\maketitle
\begin{abstract}
We describe the relation between two characterizations of conjugacy in groups of piecewise-linear
homeomorphisms, discovered by Brin and Squier in \cite{brin2} and Kassabov and Matucci in 
\cite{kama}. Thanks to the interplay between the techniques, 
we produce a simplified point of view of conjugacy that allows us to easily recover centralizers
and lends itself to generalization.
\end{abstract}

\section{Introduction}

We denote by $\PL_+(I)$ the group of orientation-preserving piecewise-linear homeomorphisms of the unit interval
$I=[0,1]$ with finitely many breakpoints. We will treat only the case of $\PL_+(I)$ even if all the results can be
adapted to certain subgroups of $\PL_+(I)$ of homeomorphisms with certain requirements on the breakpoints and the slopes
(for example, Thompson's group $F$ and the Thompson-Stein groups $\PL_{S,G}(I)$ introduced in the works
of Stein \cite{stein1} and Bieri-Strebel \cite{bieristrebel1}). In particular,
it is sufficient to restrict our study to
functions that do not intersect the diagonal, except for the points $0$ and $1$ (see Section
\ref{sec:stair-algorithm} for the motivation).

In their work \cite{brin2} Brin and Squier define an invariant under conjugacy for maps of $\PL_+(I)$ that do not intersect
the diagonal.
Their description is based on similar earlier work by Mather \cite{Math}
for diffeomorphisms of the unit interval and allows the classification of centralizers and the detection of roots
of elements. These techniques were originally introduced as an attempt to 
solve the conjugacy problem in Thompson's group $F$ (which was then proved to be solvable by Guba and Sapir in
\cite{gusa1}). Later on this approach was refined by Gill and Short in \cite{ghisho1} and Belk and Matucci
(see \cite{matucci}) to give another proof of the solution to the conjugacy problem in Thompson's group $F$.
On the other hand, Kassabov and Matucci showed a solution to the simultaneous conjugacy problem in \cite{kama}
by producing an algorithm to build all conjugators, if they exist. Similarly, these techniques can be used to
obtain centralizers and roots as a byproduct.

The aim of this note is to show the connection between the techniques in \cite{brin2} and \cite{kama}
to characterize conjugacy in groups of piecewise-linear homeomorphisms. By defining a modified version
of Brin and Squier's invariant and using a mixture of those points of view
it is possible to produce a short proof of the description of conjugacy and centralizers in $\PL_+(I)$.
In particular, the interplay between these two points of view lends itself to
generalizations giving a tool to study larger class of groups of piecewise-linear homeomorphisms.

This paper is organized as follows. In Section \ref{sec:stair-algorithm}
we give a short account of
a key algorithm in \cite{kama} (the \emph{stair algorithm}) to build a particular conjugator $g$
for two elements $y,z \in \PL_+(I)$.
In Section \ref{sec:mather-definition} we define a conjugacy invariant (called \emph{Mather invariant})
that essentially encodes the 
characterization of conjugacy in \cite{brin2} for $\PL_+(I)$. In Section \ref{sec:equivalence} we show to use the 
stair algorithm to simplify the proof of the
the characterization of conjugacy of \cite{brin2} using Mather invariants. In turn,
in Section \ref{sec:applications-generalizations} we will show how Mather invariants allow us to shorten the arguments
in \cite{kama} to classify centralizers of elements. We finish by briefly describing possible extensions of these tools.

\section{The stair algorithm for functions in $\PL_+^<(I)$ \label{sec:stair-algorithm}}

In this Section we will discuss how to find a special conjugator $g \in \PL_+(I)$ 
for two functions $y,z \in \PL_+(I)$, if it exists. The idea will be to assume that such a conjugator $g$
exists and obtain conditions that $g$ must satisfy.

\definition{We denote by $\PL_+^<(I)$ the subset of $\PL_+(I)$ of all functions that lie
below the diagonal, that is the maps $z \in \PL_+(I)$ such that $f(t)<t$ for all $t \in (0,1)$.
Similarly, we define the subset $\PL_+^>(I)$ of functions that lie above the diagonal.
A function $z \in \PL_+(I)$ is defined to be a \emph{one-bump function} if either $z \in \PL_+^<(I)$ or
$z \in \PL_+^>(I)$.}

We will restrict to study conjugacy for one-bump functions. 
The reason for this assumption is easily explained: if two functions $y,z \in \PL_+(I)$ are conjugate through
$g$, then $g^{-1}(\partial \mathrm{Fix}(y))=\partial \mathrm{Fix}(g^{-1}yg)=\partial \mathrm{Fix}(z)$; since the boundary of the set of fixed points
of either $y$ or $z$ is finite, the first step to verify conjugacy is to check if 
$\partial \mathrm{Fix}(y)$ and $\partial \mathrm{Fix}(z)$ have the same size. If this is the case,
we can always build a map $h \in \PL_+(I)$ such that $h^{-1} \left(\partial \mathrm{Fix}(y)\right)=\partial \mathrm{Fix}(z)$,
hence we reduce to check if $h^{-1}yh$ and $z$, which share the same boundary of the fixed set, are conjugate;
this is true if, for any two consecutive points $t_i,t_{i+1} \in \partial \mathrm{Fix}(z)$, we can
find a conjugator $g_i \in \PL_+([t_i,t_{i+1}])$ for the restrictions of $h^{-1}yh$ and $z$ to $[t_i,t_{i+1}]$,
which are either identity maps or one-bump functions. By restricting the study of conjugacy
to the intervals $[t_i,t_{i+1}]$, we derive our assumption on the maps.

If $z \in \PL_+(I)$, we define \emph{initial slope} and \emph{final slope}, respectively, to be the numbers $z'(0)$
and $z'(1)$. It is clear that if two one-bump functions $y$ and $z$ are conjugate, their
initial and final slope are the same. A more interesting fact is that a conjugator has to be linear in certain boxes around $0$ and
$1$. This fact, together with the ability to identify the two functions step by step, allows us to build a conjugator.

\lemmaname{Kassabov and Matucci, \cite{kama}}{Suppose $y,z \in \PL_+^<(I)$.
\begin{enumerate}
\item (\textbf{initial box}) Let $g \in \PL_+(I)$ be such that $g^{-1}yg=z$. Assume $y(t)=z(t)=ct$ for 
$t \in [0,\alpha]$ and $c<1$. Then the graph of $g$ is linear inside the box $[0,\alpha]\times[0,\alpha]$.
A similar statement is true for a ``final box''.
\item (\textbf{identification trick}) Let $\alpha \in (0,1)$ be such that $y(t)=z(t)$ for $t \in [0,\alpha]$.
Then there exists a $g \in \PL_+(I)$ such that $z(t)=g^{-1}y g(t)$ for $t \in [0,z^{-1}(\alpha)]$
and $g(t)=t$ in $[0,\alpha]$. The element $g$ is uniquely defined up to the point
$z^{-1}(\alpha)$.
\item (\textbf{uniqueness of conjugators}) For any positive real number $q$ there exists at most
one $g \in \PL_+(I)$ such that $g^{-1}yg=z$ and $g'(0)=q$.
\item (\textbf{conjugator for powers})
Let $g \in \PL_+(I)$ and $n \in \mathbb{N}$. Then $g^{-1} y g =z$ if and only if $g^{-1} y^n g =z^n$.
\end{enumerate} \label{thm:conjugacy-tools}}

\noindent \emph{Proof.} The proof of (1) is straightforward. To prove (2) we observe that, 
if such a $g$ exists then, for $t \in [0,z^{-1}(\alpha)]$
$$
y(g(t))= g(z(t))=z(t)
$$
since $z(t) \le \alpha$. Thus $g(t)=y^{-1} z(t)$ for $t \in [0,z^{-1}(\alpha)]$. 
To prove that such a $g$ exists, define
$$
g(t):=
\begin{cases}
t & t \in [0,\alpha] \\
y^{-1} z(t) & t \in [\alpha,z^{-1}(\alpha)]
\end{cases}
$$
and extend it to $I$ as a line from the point $(z^{-1}(\alpha),y^{-1}(\alpha))$ to $(1,1)$. 
To prove (3), assume that there exist two conjugators $g_1,g_2$
with initial slope $q$. Since $g_1^{-1} y g_1=g_2^{-1} y g_2$ we have that $g:=g_1 g_2^{-1}$ centralizes $y$
and it has initial slope $1$. Assume, by contradiction, that $g$ is the identity on $[0,\alpha]$ for
some $\alpha$, but $g'(\alpha^+) \ne 1$. Since we have 
\[
y(g(t))= g(y(t))=y(t)
\]
for $t \in [\alpha,y^{-1}(\alpha)]$, this implies that $g(t)=y^{-1}y(t)=t$ on $[\alpha,y^{-1}(\alpha)]$,
which is a contradiction. To prove the last statement we observe that if $f:=g^{-1} y^n g =z^n$,
then $f$ is centralized by both $g^{-1} y g$ and $z$. Since $g^{-1} y g$ and $z$ have the same
initial slope, then by (3) we have $g^{-1} y g =z$. $\square$

\bigskip
Part (1) of the previous Lemma tells us that any given conjugator $g$ must be linear in two suitable boxes
$[0,\alpha]^2$ and $[\beta,1]^2$, hence if we are given a point $(p,g(p))$ in any of those
boxes (say the final one), we can draw the longest segment contained in $[\beta,1]^2$ passing
through $(p,g(p))$ and $(1,1)$ and obtain the map $g$ in that box.
We are now going to build a candidate conjugator with a given initial slope.

\theoremname{Stair Algorithm, \cite{kama}}{Let $y,z \in \PL_+^<(I)$, let $[0,\alpha]^2$ be the initial linearity box
and let $0<q<1$ be a real number. There is an $N \in \mathbb{N}$ such that the unique candidate conjugator
with initial slope $q$ is given by
\[
g(t)=y^{-N} g_0z^N(t) \qquad \forall t \in [0,z^{-N}(\alpha)]
\]
and linear otherwise, where $g_0$ is any map in $\PL_+(I)$ which is linear in the initial box and
such that $g_0'(0)=q$. 
\label{thm:explicit}}

By ``unique candidate conjugator'' we mean a function $g$ such that, if there exists a conjugator
between $y$ and $z$ with initial slope $q$, then it must be equal to $g$. Hence we can test our candidate conjugator
to verify if it is indeed a conjugator.

\medskip
\noindent \emph{Proof.} Let $[\beta,1]^2$ be the final box and $N$ an integer big enough so that
\[
\min\{z^{-N}(\alpha),y^{-N}(q\alpha)\} > \beta.
\]
We will build a candidate conjugator $g$ between $y^N$ and $z^N$ (if it exists) 
as a product of two functions $g_0$ and $g_1$. We note that the linearity
boxes for $y^N$ and $z^N$ are still given by $[0,\alpha]^2$ and $[\beta,1]^2$.
By Lemma \ref{thm:conjugacy-tools}(1) $g$ has to be linear on $[0,\alpha]$ and so we define
an ``approximate conjugator'' $g_0$ by:
$$
g_0(t):= qt  \qquad t \in [0,\alpha]
$$
and extend it to the whole $I$ as a line through $(1,1)$. 
We then define $y_1:=g_0^{-1} y g_0$ and look for a conjugator $g_1$ of $y_1^N$ and $z^N$,
noticing that $y_1^N$ and $z^N$ coincide on $[0,\alpha]$. By the proof of Lemma \ref{thm:conjugacy-tools}(3), we define
$$
g_1(t):=
\begin{cases}
t                & \; \; t \in [0,\alpha] \\
y_1^{-N} z^N(t)       & \; \; t \in [\alpha,z^{-N}(\alpha)]
\end{cases}
$$
and extend it to $I$ as a line through $(1,1)$ so that $g_1^{-1}y_1^N g_1=z^N$ on
$[0,z^{-N}(\alpha)]$. Finally, build a function $g$ such that
$g(t):=g_0g_1(t)$ for $t \in [0,z^{-N}(\alpha)]$ and extend it to $I$ as a line through $(1,1)$
on $[z^{-N}(\alpha),1]$. The map $g$ is inside the final box at $t=z^{-N}(\alpha)>\beta$,
in fact
\[
g(z^{-N}(\alpha))=g_0 g_0^{-1} y^{-N} g_0 (\alpha)	 =   y^{-N}(q\alpha) >\beta.
\]
We observe that, by construction, $g$ is a conjugator
for $y^N$ and $z^N$ on $[0,z^{-N}(\alpha)]$, that is $g=g_0g_1=y^{-N} g_0g_1z^N$ on $[0,z^{-N}(\alpha)]$. 
Therefore
\[
g(t)=y^{-N} g_0g_1z^N(t)=y^{-N} g_0 z^N(t) \qquad \forall t \in [0,z^{-N}(\alpha)]
\]
since $g_1 z^N(t)=z^N(t)$ for $t \in [0,z^{-N}(\alpha)]$. 

By parts (1) and (3) of Lemma \ref{thm:conjugacy-tools},
if there is a conjugator for $y^N$ and $z^N$ with initial slope $q$, it must be equal to $g$. So we just
check if $g$ conjugates $y^N$ to $z^N$. Moreover, Lemma \ref{thm:conjugacy-tools}(4) tells
us that $g$ is a conjugator for $y^N$ and $z^N$ if and only it is for $y$ and $z$ and so we are done.
We remark that this proof does not depend on the choice of $g_0$. 
The only requirements on $g_0$ are that it must
be linear in the initial box and $g_0'(0)=q$.
$\square$

\section{Mather invariants for functions in $\PL_+^>(I)$ \label{sec:mather-definition}}

In this Section we will give an alternate description of Brin and Squier's conjugacy
invariant in \cite{brin2}. This reformulation was also used by Belk and Matucci (see \cite{matucci}) to
characterize conjugacy in Thompson's group $F$: however, their proof relies on special kinds
of diagrams peculiar to $F$ and cannot be generalized to other groups of homeomorphisms.

Roughly speaking, the Mather invariant of a map $z \in \PL_+^>(I)$ 
is defined by taking a power of $z$ large enough so that points very close to $0$ get mapped
to points very close to $1$. 

We will now define it precisely. Consider a one-bump function $z\in \PL_+^>(I)$, with initial slope $m_0$ and
final slope $m_1$.  In a neighborhood of zero, $z$ acts as multiplication by $m_0$:
for any sufficiently small $t > 0$ and sufficiently small powers of $z$, 
we have $z(t)=m_0 t, z^2(t)=m_0^2 t, z^3(t)=m_0^3 t, \ldots$, that is the interval
$[t,m_0t]$ is a ``fundamental domain'' for the action of $z$:
\begin{figure}[0.5\textwidth]
\includegraphics[width=6cm]{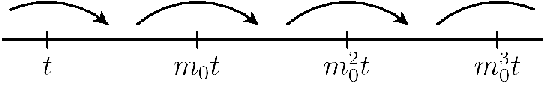}
\centering
\end{figure}

If we make the identification $t \sim m_0t$ in the interval $(0, \epsilon)$, for a sufficiently small
$\epsilon>0$, we obtain a circle $C_0$, with natural projection map
$p_0 \colon (0, \epsilon)\rightarrow C_0$. 
Similarly, if we identify $(1 - t) \sim (1 - m_1 t)$ on the interval $(1 - \delta, 1)$,
for a sufficiently small $\delta>0$, we obtain a circle $C_1$, with
natural projection map $p_1\colon (1 - \delta, 1)\rightarrow C_1$.

Let $\epsilon'>0$ be small enough so that $(\epsilon',\epsilon)$ surjects onto $C_0$:
if $N$ is sufficiently large, then $z^N$ will take $(\epsilon', \epsilon)$ 
and map it to the interval $(1 - \delta, 1)$.
This induces a map $z^\infty\colon C_0\rightarrow C_1$, making the following diagram commute:
$$
\begindc{\commdiag}[1]
\obj(-5,55){$(\epsilon', \epsilon)$}
\obj(68,55){$(1-\delta, 1)$}
\obj(0,0){$C_0$}
\obj(55,0){$C_1$}
\mor(0,55)(55,55){$z^N$}[\atleft, \solidarrow]
\mor(0,55)(0,0){$p_0$}[\atright, \solidarrow]
\mor(55,55)(55,0){$p_1$}[\atleft, \solidarrow]
\mor{$C_0$}{$C_1$}{$z^\infty$}[\atleft, \dasharrow]
\enddc
$$
The map $z^\infty$ defined above is called the \textbf{Mather invariant} for $z$.
We note that $z^\infty$ does not depend on the specific value of $N$ chosen. Any map $z^m$,
for $m \ge N$, induces the same map $z^\infty$. This is because $z$ ``acts as the identity on $C_1$'':
we can write $z^m(t)$ as $z^{m-N}(z^N(t))$, with $z^N(t) \in (1-\delta,1)$ and so, by definition of $\sim$,
we have $z^m(t) \sim z^N(t)$.
If $k > 0$, then the map $t \mapsto kt$ on $(0, \epsilon)$ induces a ``rotation'' $\rot_k$ of $C_0$. In particular, if we use the
coordinate $\theta = \log t$ on $C_0$, then
\begin{equation*}
\rot_k (\theta) \,=\, \theta + \log k
\end{equation*}
so rot$_k$ is an actual rotation. In the next Section we will give a characterization of conjugacy for one-bump functions
by means of Mather invariants.

\section{Equivalence of the two points of view \label{sec:equivalence}}

In this Section we will show the relation between the stair algorithm and the definition of Mather invariant.
This will provide an alternative proof of Brin and Squier's conjugacy invariant.

\theoremname{Brin and Squier, \cite{brin2} \label{thm:brin-mather}}{Let $y, z \in \PL_+^>(I)$ be one-bump functions with
$y'(0) = z'(0)$ and $y'(1) = z'(1)$, and let $y^\infty, z^\infty\colon C_0\rightarrow C_1$ be the corresponding Mather
invariants. Then $y$ and $z$ are conjugate if and only
if $y^\infty$ and $z^\infty$ differ by rotations of the domain and range circles \label{thm:brin-mather-invariant}:
\begin{equation*}
\begindc{\commdiag}[50]
\obj(0,1){$C_0$}
\obj(1,1){$C_1$}
\obj(0,0){$C _0$}
\obj(1,0){$C _1$}
\mor{$C_0$}{$C_1$}{$y^\infty$}[\atleft, \solidarrow]
\mor{$C _0$}{$C _1$}{$z^\infty$}[\atright, \solidarrow]
\mor{$C_0$}{$C _0$}{$\rot_k$}[\atright, \solidarrow]
\mor{$C_1$}{$C _1$}{$\rot_\ell$}[\atleft, \solidarrow]
\enddc
\end{equation*}}

\noindent \emph{Proof.} Since $y'(0)=z'(0)$ we can pick the fundamental domain for
 $y$ and $z$ around 0 to be the same. Similarly,
we can do it around 1 and so it makes sense to talk about rotations for $C_0$ and $C_1$. We stress
that the Mather invariants $y^{\infty}$ and $z^{\infty}$ that we now use
depend on the choice of the fundamental domains around 0 and 1 to talk about well defined
compositions.

We assume $z = g^{-1} y g$ for some $g \in \PL_+ (I)$
and follow the notation of the previous Section,
taking $\epsilon,\epsilon',\delta>0$ small enough and $N$ large enough. 
Then $z^N = g^{-1}y^N g$ and the following diagram commutes,
where $k = g'(0)$ and $\ell = g'(1)$:
\begin{equation*}
\begindc{\commdiag}[17]
\obj(0,0){$C_0$}
\obj(5,0){$C_1$}
\obj(2,2){$C _0$}
\obj(7,2){$C _1$}
\obj(0,5){$(\epsilon',\epsilon)$}
\obj(5,5){$(1-\delta,1)$}
\obj(2,7){$(\epsilon', \epsilon)$}
\obj(7,7){$(1-\delta, 1)$}
\mor(0,0)(5,0){$z^\infty$}[\atright, \solidarrow]
\mor(0,0)(2,2){$\rot_k$}[\atright, \solidarrow]
\mor(5,0)(7,2){$\rot_\ell$}[\atright, \solidarrow]
\mor(2,2)(5,2){$\quad y^\infty$}[\atleft, \solidline]
\mor(5,2)(7,2){}[\atleft, \solidarrow]
\mor(0,5)(0,0){$p_0$}[\atright, \solidarrow]
\mor(2,7)(2,5){}[\atleft, \solidline]
\mor(2,5)(2,2){$p_0$}[\atleft, \solidarrow]
\mor(5,5)(5,0){$\begin{matrix} p_1 \\ \, \end{matrix}$}[\atleft, \solidarrow]
\mor(7,7)(7,2){$p_1$}[\atleft, \solidarrow]
\mor{$(\epsilon',\epsilon)$}{$(1-\delta,1)$}{$\qquad z^N$}[\atleft, \solidarrow]
\mor{$(\epsilon',\epsilon)$}{$(\epsilon', \epsilon)$}{$g$}[\atleft, \solidarrow]
\mor{$(1-\delta,1)$}{$(1-\delta, 1)$}{$g$}[\atleft, \solidarrow]
\mor{$(\epsilon', \epsilon)$}{$(1-\delta, 1)$}{$y^N$}[\atleft, \solidarrow]
\enddc
\tag*{\qedhere} \; \; \;\end{equation*}

To show the converse, choose
$g_0 \in \PL_+(I)$ that is linear in the initial box and such that $g_0'(0)=k$ and define
the map $g$ to be the following pointwise limit
\[
g(t):=\lim_{n \to \infty} y^n g_0 z^{-n}(t).
\]
By the Stair Algorithm (Theorem \ref{thm:explicit}) it is clear that $g$ conjugates $y^{-1}$ to $z^{-1}$
(and hence $y$ to $z$).
It remains to show that $g \in \PL_+(I)$. By construction, $g$ has only finitely many breakpoints
on the interval $[0,1-\delta]$ for a $\delta>0$ small enough. Since $g$ conjugates $y$ and $z$,
then $g$ induces a well-defined map $g_\mathrm{ind}:C_1 \to C_1$ (given by $p_1 g p_1^{-1}$)
and we can build a diagram
similar to the one of ``only if'' part of this Theorem
\begin{equation*}
\begindc{\commdiag}[17]
\obj(0,0){$C_0$}
\obj(5,0){$C_1$}
\obj(2,2){$C _0$}
\obj(7,2){$C _1$}
\obj(0,5){$(\epsilon',\epsilon)$}
\obj(5,5){$(1-\delta,1)$}
\obj(2,7){$(\epsilon', \epsilon)$}
\obj(7,7){$(1-\delta, 1)$}
\mor(0,0)(5,0){$z^\infty$}[\atright, \solidarrow]
\mor(0,0)(2,2){$\rot_k$}[\atright, \solidarrow]
\mor(5,0)(7,2){$g_\mathrm{ind}$}[\atright, \solidarrow]
\mor(2,2)(5,2){$\quad y^\infty$}[\atleft, \solidline]
\mor(5,2)(7,2){}[\atleft, \solidarrow]
\mor(0,5)(0,0){$p_0$}[\atright, \solidarrow]
\mor(2,7)(2,5){}[\atleft, \solidline]
\mor(2,5)(2,2){$p_0$}[\atleft, \solidarrow]
\mor(5,5)(5,0){$\begin{matrix} p_1 \\ \, \end{matrix}$}[\atleft, \solidarrow]
\mor(7,7)(7,2){$p_1$}[\atleft, \solidarrow]
\mor{$(\epsilon',\epsilon)$}{$(1-\delta,1)$}{$\qquad z^N$}[\atleft, \solidarrow]
\mor{$(\epsilon',\epsilon)$}{$(\epsilon', \epsilon)$}{$g$}[\atleft, \solidarrow]
\mor{$(1-\delta,1)$}{$(1-\delta, 1)$}{$g$}[\atleft, \solidarrow]
\mor{$(\epsilon', \epsilon)$}{$(1-\delta, 1)$}{$y^N$}[\atleft, \solidarrow]
\enddc
\tag*{\qedhere} \; \; \;\end{equation*}
for suitable $\epsilon,\epsilon',\delta>0$ small enough and an integer $N$ big enough. 
By hypothesis the Mather invariants differ by
rotations of the domain and range circles, therefore we have
\[
\rot_\ell z^{\infty}= y^{\infty} \rot_k = g_\mathrm{ind} z^{\infty}
\]
and so, by cancellation, $g_\mathrm{ind}$ is a rotation by $\ell$.
To prove that $g\in \PL_+(I)$ we show that $g$ is linear around $1$ in
the following Claim:

\medskip
\noindent \emph{Claim:} 
If $g:I \to I$ is a continuous map and $p_1$ is a projection of a neighbourhood of 1 to $C_1$
such that $p_1 g p_1^{-1}$ is a well-defined map from $C_1$ to $C_1$
and it is a rotation of $C_1$,
then $g$ is linear on $(1-\delta,1]$ for a $\delta>0$ small enough.

\medskip
\noindent \emph{Proof of the Claim.} Let $\delta>0$ be small enough so that
$(1-\delta,1]$ is contained in the domain of $p_1$ and let $t \in (1-\delta,1]$. 
Following the notation from Section \ref{sec:mather-definition}, 
since  $p_1 g p_1^{-1}$ is a rotation by $\ell$, we have 
\[
g(t)=g(1-(1-t))=1-\ell m_1^r(1-t).
\]
for some integer $r$. Thus, for a $\lambda>0$ close enough to 1, we have
\[
g(1 - \lambda(1-t)) = 1 - \ell m_1^r \lambda(1-t)=
1 - \lambda(1-g(t)).
\]
By the previous equation, the function
\[
h(t):=\frac{1-g(t)}{1-t}
\]
satisfies
\[
h(t)=h(1 - \lambda(1-t))
\]
for $\lambda>0$ close enough to $1$, hence $h$
is locally constant on $(1-\delta,1]$ and therefore it is constant.
Since $h$ is constant, the map $g$ is then linear around $1$. $\square$

\remark{We have slightly abused the notation in the two cube diagrams of the previous proof: 
to simplify the exposition, we have not been careful in choosing the range sets for $g$
that still surject onto $C_0$ and $C_1$ (although it can be made precise).}

\remark{The previous proof shows that two functions $y,z$ are conjugate if and only if
the Stair Algorithm builds a linear map in the final linearity box and this
happens if and only if the two Mather invariants differ by rotations of the domain and the the range circles.
The Mather invariant thus gives the ``obstruction'' to finishing the Stair Algorithm at 1.}

\remark{We stress that the definition of Mather invariant and the construction of the stair algorithm do not really
depend upon the set of breakpoints and slopes of the maps $y$ and $z$. With little work, the two constructions and
their equivalence can be extended to Thompson-Stein groups (see also \cite{kama}).}

\section{Applications: centralizers and generalizations \label{sec:applications-generalizations}}

Given a map $f:S^1 \to S^1$, a \emph{lift} of $f$ is a map $F:\mathbb{R} \to \mathbb{R}$ such that $F(t+1)=F(t)+1$
for all $t \in \mathbb{R}$ and $F$ induces $f$ when passing the domain and the range to quotients via the
relation $\alpha \sim \alpha +1$. Given a lift, we talk about a \emph{maximal $V$-interval}
to refer to an interval $[a,b]$ such that $F$ is linear with slope $V$ on $[a,b]$ and $a,b$ are breakpoints for $F$.
We will give a short proof of the following well known result.

\theorem{Let $z \in \PL_+^>(I)$. Then the centralizer subgroup $C_{\PL_+(I)}(z)=\{g \in \PL_+(I) \mid gz=zg\}$
is isomorphic to the infinite cyclic group.}

\noindent \emph{Proof.} Define the following group homomorphism:
$$
\begin{array}{lrcl}
\varphi_z: & C_{\PL_+(I)}(z) & \longrightarrow & (\mathbb{R},+) \\
           & g              & \longmapsto     & \log g'(0).
\end{array}
$$
Lemma \ref{thm:conjugacy-tools}(3) implies that $\varphi_z$ is injective. 
By Theorem \ref{thm:brin-mather-invariant} any function $g$
centralizing $z$ induces two rotations $\rot_\ell,\rot_k$ such that
\[
\rot_\ell z^{\infty}= z^{\infty} \rot_k
\]
where $k = g'(0)$ and $\ell=g'(1)$. Observe that $R_\ell(t)=t+\log \ell$ and $R_k(t)=t+\log k$ are lifts of 
the two rotations $\rot_\ell,\rot_k$.
Choose a lift $Z:\mathbb{R} \to \mathbb{R}$ of $z^\infty$. The previous equality implies:
\[
Z(t)+\log \ell = R_\ell(Z(t))=Z(R_k(t)) = Z(t+\log k)
\]
which means that the graph of $Z$ can be shifted ``diagonally'' onto itself. 
The map $Z$ is piecewise-linear and, 
for any positive number $r$, has finitely many breakpoints on the interval $[-r,r]$. Hence $Z$
has only finitely many maximal $Z'(0)$-intervals that are contained
in $[-r,r]$ and so there is only a 
discrete set of shifts (that is, values of $\log k=\varphi_z(g)$) which maps the graph of $Z$ onto itself, 
unless $Z$ is a line. 

To see that this is not the case, we show that $z^{\infty}$
has breakpoints.
Let $N$ be a power large enough so that a fundamental domain near $0$
is sent near 1 so that $z^N$ induces $z^\infty$, then
either $z^{-N}$ or $z^{-2N}$ has a breakpoint in the final box $[\beta,1]^2$ (this implies immediately that
$z^\infty$ must have breakpoints). If 
they were both linear, by applying the chain rule on $z^{-2N}=z^{-N} \circ z^{-N}$ first at $\beta$
and then at 1, one sees that the slope
$z^{-2N}$ on $[\beta,1]$ is simultaneously equal to the product of the slopes $z'(0)^{-N} (z^{-N})'(\beta^+)$
and $z'(1)^{-N} (z^{-N})'(1)$ and this is impossible since 
$(z^{-N})'(\beta^+)=(z^{-N})'(1)$, but $z'(0)<z'(1)$.

We have thus proved that 
the image of $\varphi_z$ must be a discrete subgroup of $(\mathbb{R},+)$ and so, by a standard fact,
it is isomorphic to $\mathbb{Z}$. $\square$

\medskip
The Mather invariant approach is also interesting because it lends itself to generalizations.
Let $\PL_\mathrm{dis}(\mathbb{R})$ the group of all 
orientation-preserving piecewise-linear homeomorphisms of the real line with a discrete set
of breakpoints and let $\mathrm{EP}$ be the subgroup of $\PL_\mathrm{dis}(\mathbb{R})$ of the functions that are 
``eventually periodic at infinity'', that is functions $f \in \PL_\mathrm{dis}(\mathbb{R})$ such that there exist numbers
$L_f,R_f$ so that $f(t-1)=f(t)-1$ for $t<L_f$ and $f(t+1)=f(t)+1$ for $t>R_f$. It is easy to define the subset
$\mathrm{EP}^>$ and Mather invariant for functions in $\mathrm{EP}^>$: 
we just mod out the intervals $(-\infty,L_f)$ and 
$(R_f,\infty)$ by the relation $t \sim f(t)$ and then take a power of $f$ high enough so that $(f^{-1}(L_f),L_f)$
gets carried to a subset of $(R_f,\infty)$. Similarly, one can partially extend the stair algorithm 
to build conjugators. It is thus interesting to see how much of these techniques can
be extended to overgroups containing $\PL_+(I)$ to compute centralizers and, possibly, to
study the conjugacy problem.


\section*{Acknowledgments} 
The author would like to thank Ken Brown, Jos\'e Burillo, Martin Kassabov and an anonymous referee
for helpful comments that improved the presentation of this paper.
\bibliographystyle{plain}

\small{\noindent  \textsc{
\rule{0mm}{6mm} \\
Francesco Matucci 
\\ Centre de Recerca Matem\`atica, 
\\ Apartat 50, 08193 Bellaterra, Barcelona, Spain} 
\\ \emph{E-mail address:} \texttt{fmatucci@crm.cat}}

\end{document}